\documentclass[a4paper,12pt]{article}
\usepackage{comment}
\usepackage{cite}
\usepackage{amsmath}
\usepackage{amssymb}
\usepackage{amsfonts}
\usepackage[T1]{fontenc}
\usepackage[utf8]{inputenc}
\usepackage{graphicx}
\usepackage{fancyhdr}
\usepackage{float}
\usepackage{xcolor}
\usepackage{authblk}
\usepackage{mathrsfs}
\usepackage{empheq}
\usepackage[hyphens]{url}
\usepackage{hyperref} 
\usepackage[]{breakurl}

\usepackage{geometry}
\begin{document}

\title{Normalized centered moments of the Fr\'echet extreme-value distribution and inference of its parameter}

\author[$\dagger$]{Jean-Christophe {\sc Pain}$^{1,2,}$\footnote{jean-christophe.pain@cea.fr}\\
\small
$^1$CEA, DAM, DIF, F-91297 Arpajon, France\\
$^2$Universit\'e Paris-Saclay, CEA, Laboratoire Mati\`ere en Conditions Extr\^emes,\\ 
91680 Bruy\`eres-le-Ch\^atel, France
}

\maketitle

\begin{abstract}
In the present work, we provide the general expression of the normalized centered moments of the Fr\'echet extreme-value distribution. In order to try to represent a set of data corresponding to rare events by a Fréchet distribution, it is important to be able to determine its characteristic parameter $\alpha$. Such a parameter can be deduced from the variance (proportional to the square of the Full Width at Half Maximum) of the studied distribution. However, the corresponding equation requires a numerical resolution. We propose two simple estimates of $\alpha$ from the knowledge of the variance, based on the Laurent series of the Gamma function. The most accurate expression involves the Ap\'ery constant.
\end{abstract}

\section{Introduction}

The Fr\'echet law \cite{Frechet1927} is, along with the Gumbel \cite{Gumbel1958} and Weibull \cite{Papoulis2002} ones, a type of extreme-value distribution, used to model the distribution of the maximum value of a sample. It finds many applications in different fields such as natural calamities, horse racing, rainfall, queues in supermarkets, sea currents or wind speeds (see for instance Refs. \cite{Kotz2000,Coles2001}).

The generalized Fr\'echet law (probability distribution function) reads
\begin{equation}\label{gen}
g(x;\alpha,s,m)=\frac{\alpha}{s}\left(\frac{x-m}{s}\right)^{-1-\alpha}e^{-(\frac{x-m}{s})^{-\alpha}},
\end{equation}
depending on parameters $m$ (position of the maximum), $s>0$ (scale parameter) and $\alpha$ the shape parameter. The associated cumulative distribution function is $e^{-\left(\frac{x-m}{s}\right)^{-\alpha}}$ if $x>m$ and 0 otherwise. For simplicity, in the present work we focus on the one-parameter Fr\'echet distribution
\begin{equation}
f(x;\alpha)=\alpha~x^{-\alpha-1}e^{-x^{-\alpha}}, 
\end{equation}
with $\alpha>0$, corresponding to the cumulative distribution function $e^{-x^{\alpha}}$. The results presented here can be easily generalized to the form of Eq. (\ref{gen}).

It can be useful to fit an observed (or measured) distribution of rare events, by a Fr\'echet distribution \cite{Ramos2019,Pain2023}. To do so, we need to determine the parameter $\alpha$. The expression of the normalized centered moments of the Fr\'echet distribution are  given in section \ref{sec2}. In section \ref{sec3}, we provide an approximate determination of the parameter $\alpha$ from the knowledge of the variance of a distribution.

\section{Normalized centered moments}\label{sec2}

The moments of the one-parameter Fr\'echet distribution are respectively (see for instance Ref. \cite{Muraleedharan2009}):
\begin{equation}
\mu_{k}=\int_{0}^{\infty}x^{k}f(x;\alpha)\,dx=\int_{0}^{\infty }t^{-{\frac {k}{\alpha }}}e^{-t}\,dt=\Gamma \left(1-{\frac {k}{\alpha}}\right),
\end{equation}
where $\Gamma$ is the usual Gamma function and $k\geq 1$. The moments $\mu_k$ are defined for $k<\alpha$. The centered moments are defined as
\begin{equation}
\mu_{k,c}=\int_{0}^{\infty}(x-\mu_1)^{k}f(x;\alpha)\,dx.
\end{equation}
Using the binomial expansion theorem, one gets
\begin{equation}
\mu_{k,c}=\sum_{p=0}^k\binom{k}{p}\int_{0}^{\infty}x^k(-\mu_1)^{k-p}f(x;\alpha)\,dx,
\end{equation}
where $\binom{k}{p}=k!/p!/(k-p)!$ is the usual binomial coefficient, \emph{i.e.}
\begin{equation}
\mu_{k,c}=\sum_{p=0}^k\binom{k}{p}\int_{0}^{\infty}x^k(-1)^{k-p}\left[\Gamma\left(1-\frac{1}{\alpha}\right)\right]^{k-p}f(x;\alpha)\,dx.
\end{equation}
Setting $\Omega_k=\Gamma\left(1-\frac{k}{\alpha}\right)$, the latter equation reads
\begin{equation}
\mu_{k,c}=(-1)^{k}\sum_{p=0}^k(-1)^{p}\binom{k}{p}\Omega_1^{k-p}\int_{0}^{\infty}x^kf(x;\alpha)\,dx.
\end{equation}
The reduced centered moments are defined (for $k\geq 2$) as
\begin{equation}
\zeta_{k,c}=\frac{\mu_{k,c}}{(\mu_{2,c})^{k/2}},
\end{equation}
and thus
\begin{empheq}[box=\fbox]{align}
\zeta_{k,c}=\frac{\displaystyle(-1)^{k}\sum_{p=0}^k(-1)^{p}\binom{k}{p}\Omega_1^{k-p}\int_{0}^{\infty}x^kf(x;\alpha)\,dx}{\displaystyle\left(\Omega_2-\Omega_1^2\right)^{n/2}}.
\end{empheq}
The skewness of the Fr\'echet distribution (characterizing its asymmetry) is
\begin{equation}
\zeta_{3,c}=\frac{\Omega_3-3\Omega_2\Omega_1+2\Omega_1^3}{\sqrt{\Omega_2-\Omega_1^2}^{3}},
\end{equation}
for $\alpha>3$ and $+\infty$ otherwise. The excess kurtosis (kurtosis minus three), characterizing the sharpness of the distribution as compared to the Gaussian, reads
\begin{equation}
\zeta_{4,c}-3=-6+\frac{\Omega_4-4\Omega_3\Omega_1+3\Omega_2^2}{\Omega_2-\Omega_1^2}
\end{equation}
for $\alpha>4$ and $\infty$ otherwise. In addition, for instance, the normalized centered sixth-order moment reads 
\begin{equation}
\zeta_{6,c}={\frac {\Omega_6-6\Omega_5\Omega_1+15\Omega_4\Omega_2^2-20\Omega_3\Omega_2^3+15\Omega_2\Omega_1^4-5\Omega_1^6}{\left(\Omega_2-\Omega_1^2\right)^{3}}}.
\end{equation}

\section{Inference of the Fr\'echet parameter from the knowledge of its variance}\label{sec3}

The purpose of this section is to provide an approximate determination of the parameter $\alpha$ of the Fr\'echet distribution from the knowledge of the variance.

\subsection{Laurent series of the Gamma function}

The Gamma function is
\begin{equation}
\Gamma(z)=\int_0^{\infty}e^{-t}t^{z-1}\,dt.
\end{equation}
Integrating by parts, one gets
\begin{eqnarray}
\Gamma(z)&=&\left.e^{-t}\frac{t^z}{z}\right|_0^{\infty}+\int_0^{\infty}e^{-t}\frac{t^z}{z}\,dt\nonumber\\
&=&\frac{1}{z}\int_0^{\infty}e^{-t}t^z\,dt=\frac{1}{z}\int_0^{\infty}e^{-t}e^{z\ln t}\,dt.
\end{eqnarray}
Expanding $e^{z\ln t}$ in Taylor series yields
\begin{equation}
\Gamma(z)=\frac{1}{z}\int_0^{\infty}e^{-t}\sum_{n=0}^{\infty}\frac{(z\ln t)^n}{n!}\,dt
\end{equation}
\emph{i.e.}
\begin{equation}
\Gamma(z)=\frac{1}{z}\sum_{n=0}^{\infty}\frac{z^n}{n!}\int_0^{\infty}e^{-t}(\ln t)^n\,dt,
\end{equation}
which gives 
\begin{equation}\label{expan}
\Gamma(z)=\frac{1}{z}\int_0^{\infty}e^{-t}\ln t\,dt+\frac{z}{2}\int_0^{\infty}e^{-t}(\ln t)^2\,dt+\frac{z^2}{6}\int_0^{\infty}e^{-t}(\ln t)^3\,dt+O(z^3).
\end{equation}
One knows that
\begin{equation}
\int_0^{\infty}e^{-t}\ln t\,dt=\Gamma'(1)=-\gamma.
\end{equation}
One has $\Gamma'(x)=\Gamma(x)~\psi(x)$, where $\psi$ represents the Digamma function
\begin{equation}\label{dig}
{\displaystyle \psi (z+1)=-\gamma +\sum _{n=1}^{\infty }{\frac {z}{n(n+z)}}\qquad z\neq -1,-2,-3,\ldots }
\end{equation}
From Eq. (\ref{dig}), we deduce $\psi(1)=-\gamma$, $\psi'(1)=\zeta(2)=\pi^2/6$, where $\zeta$ is the usual zeta function
\begin{equation}
\zeta(s)=\sum_{n=1}^\infty{\frac1{n^s}} = 1 + \frac1{2^s} + \frac1{3^s} + \frac1{4^s} + \cdots.
\end{equation}
The second integral in Eq. (\ref{expan}) is equal to
\begin{equation}
\int_0^{\infty}e^{-t}(\ln t)^2\,dt=\Gamma''(1)=\Gamma'(1)\psi(1)+\Gamma(1)\psi'(1)=\gamma^2+\frac{\pi^2}{6}
\end{equation}
and since $\psi''(1)=-2~\zeta(3)$, $\zeta(3)$ being the Ap\'ery constant, we have also
\begin{eqnarray}
\int_0^{\infty}e^{-t}(\ln t)^3\,dt&=&\Gamma'''(1)=\Gamma''(1)\psi(1)+2\Gamma'(1)\psi'(1)+\Gamma(1)\psi''(1)\nonumber\\
&=&-\gamma\left(\gamma^2+\frac{\pi^2}{6}\right)-2\gamma\frac{\pi^2}{6}-2~\zeta(3).
\end{eqnarray}
Finally, the Laurent expansion of the Gamma function up to order two reads
\begin{equation}\label{Laurent}
\Gamma(z)=\frac{1}{z}-\gamma+\frac{1}{2}\left(\gamma^2+\frac{\pi^2}{6}\right)z-\frac{1}{6}\left(\gamma^3+\frac{\gamma\pi^2}{2}+2~\zeta(3)\right)z^2+O(z^3).
\end{equation}

\subsection{Simple approximate formula for the parameter of the Fr\'echet distribution}

Truncating the Laurent series above (Eq. (\ref{Laurent})) at order $z$:
\begin{equation}
\Gamma(z)=\frac{1}{z}-\gamma+\frac{1}{2}\left(\gamma^2+\frac{\pi^2}{6}\right)z+O(z^2)
\end{equation}
gives, using $\Gamma(z+1)=z~\Gamma(z)$:
\begin{equation}
\Gamma(z+1)=1-\gamma~z+\frac{1}{2}\left(\gamma^2+\frac{\pi^2}{6}\right)z^2+O(z^3).
\end{equation}
Under that assumption, the variance becomes
\begin{equation}
\mathscr{V}=\Omega_2-\Omega_1^2\approx 1+\frac{2\gamma}{\alpha}+\frac{1}{2}\left(\gamma^2+\frac{\pi^2}{6}\right)\frac{4}{\alpha^2}-\left(1+\frac{\gamma}{\alpha}+\frac{1}{2}\left(\gamma^2+\frac{\pi^2}{6}\right)\frac{1}{\alpha^2}\right)^2,
\end{equation}
\emph{i.e.}
\begin{equation}
\mathscr{V}\approx\frac{\pi^2}{6\alpha^2}.
\end{equation}
The parameter $\alpha$ is thus simply
\begin{empheq}[box=\fbox]{align}\label{approx1}
\alpha\approx\frac{\pi}{\sqrt{6~\mathscr{V}}}.
\end{empheq}

\begin{table}[ht]
\centering
\begin{tabular}{|c|c|c|}\hline
Variance & $\alpha$ [exact] & $\alpha$ [formula (\ref{approx1})]\\\hline\hline
$\Gamma\left(1-\displaystyle\frac{2}{5}\right)-\Gamma^2\left(1-\displaystyle\frac{1}{5}\right) \approx 0.133761$ & 5 & 3.51\\\hline
$\Gamma\left(1-\displaystyle\frac{2}{10}\right)-\Gamma^2\left(1-\displaystyle\frac{1}{10}\right) \approx 0.0222624$ & 10 & 8.60\\\hline
$\Gamma\left(1-\displaystyle\frac{2}{50}\right)-\Gamma^2\left(1-\displaystyle\frac{1}{50}\right) \approx 0.000694362$ & 50 & 48.67\\\hline
$\Gamma\left(1-\displaystyle\frac{2}{100}\right)-\Gamma^2\left(1-\displaystyle\frac{1}{100}\right) \approx 0.000168916$ & 100 & 98.68\\\hline\hline
\end{tabular}
\caption{Comparison between the exact value of $\alpha$ and estimate (\ref{approx1}).}\label{tab1}
\end{table}

The accuracy of expression (\ref{approx1}) can be appreciated from table \ref{tab1}.

\subsection{Improving the estimate using the Laurent series of $\Gamma(z)$ up to $z^2$}
Using the Laurent series (\ref{Laurent}), one has
\begin{equation}
\Gamma(z+1)=1-\gamma~z+\frac{1}{2}\left(\gamma^2+\frac{\pi^2}{6}\right)z^2-\frac{1}{6}\left(\gamma^3+\frac{\gamma\pi^2}{2}+2~\zeta(3)\right)z^3+O(z^4)
\end{equation}
and finally
\begin{empheq}[box=\fbox]{align}\label{approx2}
\frac{\pi^2}{6}\left(\frac{1}{\alpha}\right)^2+\frac{(\gamma\pi^2+6~\zeta(3))}{3}\left(\frac{1}{\alpha}\right)^3\approx\mathscr{V}.
\end{empheq}
The latter equation is a third-order polynomial equation in the variable $1/\alpha$; it can be solved analytically using the Tartaglia-Cardano formulas \cite{McKelvey1984}.

\begin{table}[!ht]
\centering
\begin{tabular}{|c|c|c|}\hline
Variance & $\alpha$ [exact] & $\alpha$ [solution of Eq. (\ref{approx2})]\\\hline\hline
$\Gamma\left(1-\displaystyle\frac{2}{5}\right)-\Gamma^2\left(1-\displaystyle\frac{1}{5}\right) \approx 0.133761$ & 5 & 4.42\\\hline
$\Gamma\left(1-\displaystyle\frac{2}{10}\right)-\Gamma^2\left(1-\displaystyle\frac{1}{10}\right) \approx 0.0222624$ & 10 & 9.69\\\hline
$\Gamma\left(1-\displaystyle\frac{2}{50}\right)-\Gamma^2\left(1-\displaystyle\frac{1}{50}\right) \approx 0.000694362$ & 50 & 49.93\\\hline
$\Gamma\left(1-\displaystyle\frac{2}{100}\right)-\Gamma^2\left(1-\displaystyle\frac{1}{100}\right) \approx 0.000168916$ & 100 & 99.965\\\hline\hline
\end{tabular}
\caption{Comparison between the exact value of $\alpha$ and solution of approximate equation (\ref{approx2}).}\label{tab2}
\end{table}

Table \ref{tab2} shows some estimates of $\alpha$ given by the analytical solution of Eq. (\ref{approx2}), compared to the exact values.

%\clearpage

\section{Conclusion}

We provided the general expressions of the normalized centered moments of the Fr\'echet extreme-value distribution. Representing a set of data corresponding to rare events by a Fréchet distribution requires determining its characteristic parameter $\alpha$. The latter can be deduced from the variance of the studied distribution. However, the corresponding equation needs to be solved numerically. We propose two simple estimates of $\alpha$, based on the Laurent series of the Gamma function, up to monomials proportional to $z$ and $z^2$ respectively.

\end{document}